\documentclass[a4paper,12pt]{amsart}
\usepackage{amsmath,latexsym,amscd} 
\title{Canonical stability of 3-folds of general type with $p_g\geq 3$}
\author{Meng Chen}
\thanks{The paper is partially supported by the National Natural Science Foundation of China 
(Key Project No. 10131010),  Shanghai Scientific $\&$ Technical Commission 
(Grant 01QA14042) and SRF for ROCS, SEM}
\address{Department of Applied Mathematics, Tongji University, 
Shanghai, 200092, PR China}
\email{mchen@tongji.edu.cn}

\newcommand{\roundup}[1]{\ulcorner{#1}\urcorner}
\newcommand{\rounddown}[1]{\llcorner{#1}\lrcorner}

\newtheorem{thm}{Theorem}[section]
\newtheorem{lem}[thm]{Lemma}

\newtheorem{prop}[thm]{Proposition}
\newtheorem{claim}[thm]{Claim}
\theoremstyle{definition}

\newtheorem{setup}[thm]{}
\newtheorem{question}[thm]{Question}
\newtheorem{exmp}[thm]{Example}

\newtheorem{rem}[thm]{Remark}           
\theoremstyle{remark}

\begin{document}
\begin{abstract}
We study the canonical stability of a smooth projective 3-fold $V$ of general type. 
We prove that $|5K_V|$ gives a birational map onto its image provided $p_g\ge 4$. 
For those $V$ with $p_g=3$, $|6K_V|$ gives a birational map.  
Known examples show that both the results are optimal.
\end{abstract}
\maketitle
\pagestyle{myheadings}
\markboth{\hfill M. Chen\hfill}{\hfill Canonical stability of 3-folds of general 
type\hfill}
\section{\bf Introduction}
Throughout the base field is ${\Bbb C}$. Let $V$ be a smooth projective variety with 
$\kappa(V)>0$. For all integer $m>0$, one may define the so-called m-canonical map 
$\Phi_m$, which is nothing but the rational map corresponding to the complete linear 
system $|mK_V|$. Within birational geometry, to study the behavior of $\Phi_m$ has been one 
of the classical and important aspects. When $\dim V\leq 2$, the behavior of 
$\Phi_m$ is known quite well according to works by a long list of authors. When 
$\dim V\ge 3$, known results on $\Phi_m$ are partial and there remains a lot of open problems on this topic. In this paper, we study a 3-fold of general type.     

Assume $\dim V=\kappa(V)=3$. According to Mori's MMP, $V$ has a 
minimal model $X$ with only ${\Bbb Q}$-factorial terminal singularities. Denote by 
$K_X$ the canonical Weil divisor on $X$. One may also define the m-canonical map 
$\varphi_m$ corresponding to the complete linear system $|mK_X|$. Modulo birational 
equivalence, both $\Phi_m$ and $\varphi_m$ share lots of common properties, e.g. 
birationality. In general, it doesn't make sense to study the basepoint freeness of 
$|mK_X|$ unless $m$ is divisible by the canonical index $r(X)$. However, 
one may consider the birationality of $\Phi_m$ or $\varphi_m$.
It's well-known that, for a given $V$, 
$\Phi_{m(V)}$ is birational onto its image whenever 
$m(V)\gg 0$. Therefore quite an interesting thing to do is to answer the following question which is still open

\begin{question}\label{Q1} Let $V$ be a smooth projective 3-fold of general type. 

(1) Do there exist a universal lower bound $N$ ($N$ doesn't depend on $V$) such that $\Phi_m$ is birational onto its image for all $m\ge N$?

(2) For a given $V$, what is the optimal lower bound $m_0(V)$ such that $\Phi_{m_0(V)}$ 
is birational?
\end{question}   

Apart from the generally accepted importance of 1.1(1) above, we would like to 
emphasize here that 1.1(2) is equally important. For example, it is strongly related 
to 3-fold geography (see \cite{Hunt}) and, also, it can be applied to determine the 
automorphism group of $V$ (see the Remark in \cite{X}). There has been partial 
results on Question \ref{Q1}. According to \cite{Ch1, E-L, Lee}, $\varphi_m$ is a birational morphism onto its image if $X$ is a minimal Gorenstein 3-fold of general type. For a general 3-fold $V$ with  
$q(V)>0$ or $p_k(V)\ge 2$, J. Koll\'ar (\cite{Kol}) gave a positive answer to 
\ref{Q1}(1). For regular 3-folds of general type, T. Luo (\cite{Luo}) partially 
answered \ref{Q1}(1). There is, however, hardly no result on \ref{Q1}(2). The aim of 
this paper is to study  \ref{Q1}(2) under an extra assumption. My main result is the following

\begin{thm}\label{Main} Let $V$ be a smooth projective 3-fold of general type. 

(1) Assume $p_g(V)= 3$. Then $\Phi_6$ is birational onto its image. 

(2) Assume $p_g(V)\ge 4$. Then $\Phi_5$  is birational onto its image.    
\end{thm}

The following examples show that Theorem \ref{Main} is optimal.

\begin{exmp}\label{E1} In \cite{C-G}, Chiaruttini and Gattazzo has found a smooth projective 3-fold $V$ of general type with $q(V)=h^2({\mathcal O}_V)=0$, $p_g(V)=3$. They verified, on $V$, that $\Phi_m$ is birational if and only if $m\ge 6$ and that $\Phi_m$ is generically finite for $2\le m\le 5$. 
It is, in fact, not difficult to see that $V$ is canonically fibred by curves of genus two. 
\end{exmp}

\begin{exmp}\label{E2} Denote by $S$ a smooth minimal surface of general type with 
$(K_S^2, p_g(S))=(1,2)$. Pick up a smooth curve $C_k$ of genus $k\ge 2$. Set 
$X_k:=S\times C_k$. Then $p_g(X_k)=2k\ge 4$. It is obvious that the $\Phi_4$ of $X_k$ is not birational. This is of course a trivial example. 
\end{exmp} 

\begin{exmp} On ${\Bbb P}_{\Bbb C}^3$, take a smooth
hypersurface $S$ of degree 10. $S\sim 10H$ where $H$ is a hyperplane. Let $X$ be a double cover over
${\Bbb P}^3$ with branch locus along $S$. Then
$X$ is a nonsingular canonical model. 
$K_X^3=2$ and $p_g(X)=4$ and
$\Phi_{1}$ is a finite morphism onto ${\Bbb P}^3$ of degree $2$.
One may easily check that $\Phi_{4}$ is also a finite morphism of degree $2$. 
\end{exmp}

\section{\bf The general method}
\begin{setup}\label{S:2.1} {\bf Brief review on curves.} We recall several facts which will be applied in the paper. 

(2.1.1) Let $C$ be a smooth curve of genus $\ge 2$. Then $\varphi_m$ is an embedding for all $m\ge 3$.  

(2.1.2) Let $C$ be a smooth curve of genus $\ge 2$. Assume $D$ is a divisor on $C$. Then $K_C+D$ is very ample whenever $\deg(D)\ge 3$ and $|K_C+D|$ is basepoint free whenever $\deg(D)\ge 2$.

(2.1.3) Let $C$ be a smooth non-hyperelliptic curve. Assume $D$ is a divisor on $C$ 
with $\deg(D)\ge 2$. Then the rational map corresponding to $|K_C+D|$ gives a birational morphism onto its image. 
\end{setup}

\begin{setup}\label{S:2.2} {\bf Brief review of relevant results on surfaces.} 
We only recall those cited in the paper. Let $S$ be a smooth minimal surface of general 
type. By \cite{BPV, Bom, Catan, Ci}, one has

(2.2.1) $\varphi_m$ is a birational morphism onto its image for all $m\ge 5$.

(2.2.2) $\varphi_4$ is a birational morphism if and only if $(K_S^2, p_g(S))\neq (1,2).$

(2.2.3) $\varphi_3$ is birational if and only if $(K_S^2, p_g(S))\ne (1,2)$ and $(2,3)$.

(2.2.4) $|2K_S|$ is basepoint free if $p_g(S)>0$.
\end{setup}

\begin{prop}\label{P:2.3}
Let $S$ be a smooth projective surface of general type with $(K_{S_0}^2, p_g(S))=(1,2)$, where 
$\sigma: S\rightarrow S_0$ is the contraction onto the minimal model. Assume $L$ is a nef and big
${\Bbb Q}$-divisor on $S$. Then the rational map corresponding to 
$|K_S+3\sigma^*(K_{S_0})+\roundup{L}|$ is birational onto its image.
\end{prop}
\begin{proof}
First, it is easy to reduce to the case that the movable part of $|\sigma^*(K_{S_0})|$ is 
basepoint free. So we may assume, from now on, that the movable part of $|\sigma^*(K_{S_0})|$ is 
basepoint free. Denote by $|G|$ the movable part of $|\sigma^*(K_{S_0})|$. Then $|G|$ is 
composed of a rational pencil of curves of genus 2 and $h^0(S, G)=2$ (see \cite{BPV}).
 Let $C\in |G|$ be a general member. 

The Kawamata-Viehweg vanishing theorem (see 2.4 below) gives the surjective map
$$H^0(S, K_S+\sigma^*(K_{S_0})+\roundup{L}+G)\longrightarrow 
H^0(C, K_C+D_0),$$
where $\deg(D_0)>0$. 
This means that $K_S+\sigma^*(K_{S_0})+\roundup{L}+G\ge 0$ and that 
$K_S+3\sigma^*(K_{S_0})+\roundup{L}\ge G$. So 
$|K_S+3\sigma^*(K_{S_0})+\roundup{L}|$ separates different general members of $|G|$. 

Applying the vanishing theorem again, we get the surjective map
$$H^0(S, K_S+2\sigma^*(K_{S_0})+\roundup{L}+G)\longrightarrow 
H^0(C, K_C+D),$$
where $deg(D)\ge 3$. So  $|K_S+3\sigma^*(K_{S_0})+\roundup{L}||_C$ gives a birational map.
We are done.
\end{proof}

\begin{setup}\label{Vanishing}{\bf Vanishing theorem.} We always apply the Kawamata-Viehweg 
vanishing theorem (see \cite{Ka} or \cite{V1}), which plays very effective roles throughout 
the whole context. It is also well-known that, on surfaces, one may apply the vanishing theorem
without the assumption for  "normal crossings" (see \cite{Sakai}). 
\end{setup}

\begin{setup}\label{Notations} {\bf Set up for $\varphi_1$.} Because of the 
successful 3-dimensional MMP, one may always study a minimal 3-fold. 
Let $X$ be a minimal projective 3-fold of general type with only ${\Bbb Q}$-factorial terminal singularities. 
Suppose $p_g(X)\ge 2$. We study the canonical map $\varphi_{1}$ which is usually a rational map.
 Take the birational modification $\pi: X'\rightarrow X$, according to Hironaka, such that

(i) $X'$ is smooth;

(ii) the movable part of $|K_{X'}|$ is basepoint free. (Sometimes we even call for 
such a modification that those movable parts of a finite number of linear systems are 
all basepoint free.) 

(iii) $\pi^*(K_X)$ is linearly equivalent to a divisor supported by a divisor of normal crossings.

Denote by $g$ the composition $\varphi_{1}\circ\pi$. So
$g: X'\longrightarrow W'\subseteq{\Bbb P}^{p_g(X)-1}$
is a morphism. Let
$g: X'\overset{f}\longrightarrow B\overset{s}\longrightarrow W'$
be the Stein factorization of $g$. We may write
$$K_{X'}=_{\Bbb Q}\pi^*(K_X)+E_1=_{\Bbb Q}M_1+Z_1,$$
where $M_1$ is the movable part of $|K_{X'}|$, $Z_1$ the fixed part and $E_1$ an effective ${\Bbb Q}$-divisor which is a ${\Bbb Q}$-sum of distinct exceptional divisors. Throughout we always mean $\pi^*(K_X)$ by $K_{X'}-
E_1$. We fix a divisor $K_{X'}$ from the beginning.
So, whenever we take the round up of $m\pi^*(K_X)$, we always have 
$\roundup{m\pi^*(K_X)}\le mK_{X'}$ for all positive number $m$. 
We may also write
$$\pi^*(K_X)=_{\Bbb Q} M_1+E_1',$$
where $E_1'=Z_1-E_1$ is actually an effective ${\Bbb Q}$-divisor. 

If $\dim\varphi_{1}(X)=2$, we see that a general fiber of $f$ is a smooth projective 
curve of genus $g\ge 2$. We say that $X$ is {\it canonically fibred by curves of 
genus $g$}.

If $\dim\varphi_{1}(X)=1$, we see that a general fiber $S$ of $f$ is a smooth projective surface of general type. We say that $X$ is {\it canonically fibred by surfaces} with invariants $(c_1^2(S_0), p_g(S)),$
 where $S_0$ is the minimal model of $S$.  
We may write $M_1\equiv a_1S$ where $a_1\ge p_g(X)-1$.

{\it A generic irreducible element $S$ of} $|M_1|$ means either a general member 
of $|M_1|$ whenever 
$\dim\varphi_{1}(X)\ge 2$ or, otherwise, a general fiber of $f$.  

\end{setup}

\begin{thm}\label{Key} Let $X$ be a minimal projective 3-fold of general type with only 
${\Bbb Q}$-factorial terminal singularities and assume $p_g(X)\ge 2$. Keep the same notations 
as in 2.5. Pick up a generic irreducible element $S$ of $|M_1|$. Suppose, on the smooth surface 
$S$, there is a movable linear system $|G|$ and denote by $C$ a generic irreducible element of 
$|G|$. Set $\xi:=(\pi^*(K_X)\cdot C)_{X'}$ and 
$$p:=\begin{cases} 1  &\text{if}\ \dim\varphi_1(X)\ge 2\\
a_1\geq p_g(X)-1  &\text{otherwise.}
\end{cases}$$   

Assume
    
(i) there is a positive integer $m$ such that the linear system 
$$|K_S +\roundup{(m-2)\pi^*(K_X)|_S}|$$ 
separates different generic irreducible elements of $|G|$;

(ii) there is a rational number $\beta>0$ such that $\pi^*(K_X)|_S-\beta C$ is 
numerically equivalent to an effective ${\Bbb Q}$-divisor;

(iii) either the inequality 
$\alpha:=(m-1-\frac{1}{p}-\frac{1}{\beta})\xi>1$ holds (set $\alpha_0:=\roundup{\alpha}$)
 or $C$ is non-hyperelliptic, $m-1-\frac{1}{p}-\frac{1}{\beta}>0$ and $C$ is an even divisor on $S$.

Then we have the inequality $m\xi\ge 2g(C)-2+\alpha_0$. Furthermore, 
$\varphi_m$ of $X$ is birational onto its image provided either $\alpha>2$ or 
$\alpha_0=2$  and $C$ is non-hyperelliptic or  
$C$ is non-hyperelliptic, $m-1-\frac{1}{p}-\frac{1}{\beta}>0$ and $C$ is an even divisor on $S$.
\end{thm}
\begin{proof}
We consider the sub-system 
$$|K_{X'}+\roundup{(m-1)\pi^*(K_X)-\frac{1}{p}E_1'}|\subset |mK_{X'}|.$$
This system obviously separates different generic irreducible elements of $|M_1|$. By 
the birationality principle (P1) and (P2) of \cite{Ch2}, it is sufficient to prove that $|mK_{X'}||_S$ gives a birational
map. Noting that $(m-1)\pi^*(K_X)-\frac{1}{p}E_1'-S$ is nef and big, 
the vanishing theorem gives the surjective map
\begin{eqnarray*}
&H^0(X',K_{X'}+\roundup{(m-1)\pi^*(K_X)-\frac{1}{p}E_1'})\\
\longrightarrow  &H^0(S, K_S+\roundup{(m-1)\pi^*(K_X)-S-\frac{1}{p}E_1'}|_S).
\end{eqnarray*}
We are reduced to prove that $|K_S+\roundup{(m-1)\pi^*(K_X)-S-\frac{1}{p}E_1'}|_S|$ gives a birational
map. We still apply the principle (P1) and (P2) of \cite{Ch2}. Because
$$ K_S+\roundup{(m-1)\pi^*(K_X)-S-\frac{1}{p}E_1'}|_S\ge 
K_S+\roundup{(m-2)\pi^*(K_X)|_S},$$
the linear system  $|K_S+\roundup{(m-1)\pi^*(K_X)-S-\frac{1}{p}E_1'}|_S|$
separates different irreducible elements of $|G|$ by assumption (i).  
Now pick up a generic irreducible element $C\in |G|$. By assumption (ii), there is an effective 
${\Bbb Q}$-divisor $H$ on $S$ such that
$$\frac{1}{\beta}\pi^*(K_X)|_S\equiv C+H.$$ 
By the vanishing theorem, we have the surjective map 
$$H^0(S, K_S+\roundup{(m-1)\pi^*(K_X)-S-\frac{1}{p}E_1'}|_S-H)\longrightarrow 
H^0(C, D),$$
where $D:=\roundup{((m-1)\pi^*(K_X)-S-\frac{1}{p}E_1')|_S-C-H}|_C$ is a divisor on $C$.
Noting that
$$((m-1)\pi^*(K_X)-S-\frac{1}{p}E_1')|_S-C-H\equiv (m-1-\frac{1}{p}-\frac{1}{\beta})\pi^*(K_X)|_S$$
and that $C$ is nef on $S$, we have 
$\deg(D)\geq\alpha$ and thus $\deg(D)\geq \alpha_0$.
Whenever $C$ is non-hyperelliptic, $m-1-\frac{1}{p}-\frac{1}{\beta}>0$ and $C$ is an even divisor on $S$,
 $\deg(D)\geq 2$ automatically follows.
If $\deg(D)\ge 3$, then $|K_S+\roundup{((m-1)\pi^*(K_X)-S-\frac{1}{p}E_1')|_S-H}||_C$ gives a 
birational map by (2.1.2). Because 
\begin{eqnarray*}
&|K_S+\roundup{((m-1)\pi^*(K_X)-S-\frac{1}{p}E_1')|_S-H}|\\
\subset &|K_S+\roundup{(m-1)\pi^*(K_X)-S-\frac{1}{p}E_1'}|_S|,
\end{eqnarray*}
we see that the right linear system in above gives a birational map. So $\varphi_m$ of $X$ is birational.

Whenever $\deg(D)\ge 2$, $|K_C+D|$ is basepoint free. Denote by $|M_m|$ the movable part of 
$|mK_{X'}|$ and by $|N_m|$ the movable part of 
$|K_S+\roundup{((m-1)\pi^*(K_X)-S-\frac{1}{p}E_1')|_S-H}|$. By Lemma 2.7 of \cite{Ch2}, we have
$$m\pi^*(K_X)|_S\ge N_m\ \ \text{and}\ \ (N_m\cdot C)_S\ge 2g(C)-2+\alpha_0.$$
The theorem is proved.
\end{proof}

\section{\bf Proof of the main theorem}
We study $\varphi_m$ according to the value $d:=\dim\varphi_1(X)$. Obviously $1\le d\le 3.$

\begin{thm}\label{T:3.1} Let $X$ be a minimal 3-fold of general type with only 
${\Bbb Q}$-factorial terminal singularities. Assume $p_g(X)\ge 3$ and $d\ge 2$. Then 

(i) $\varphi_5$ is birational onto its image for $p_g(X)\ge 4$.

(ii) $\varphi_6$ is birational onto its image for $p_g(X)=3$.
\end{thm}
\begin{proof} 
We take $m=5$ and will apply Theorem \ref{Key}. In this case, $S$ is a general member of $|M_1|$.
So $S$ is a smooth surface of general type. Set $G:=M_1|_S$. Then $|G|$ is a basepoint free
system on $S$, since $g$ is a morphism. Take a generic irreducible element $C$ of $|G|$. $C$ is
a smooth curve of genus $\ge 2$. We have $p=1$. 

Case 1. If $d=3$, then we may take $\beta=1$. Because $\pi^*(K_X)|_S\ge G$ and $|G|$ is not composed of a pencil
of curves, Theorem \ref{Key}(i) is satisfied.  Noting that $C^2\ge 2$ in this case, 
we have 
$$\xi=(\pi^*(K_X)|_S\cdot C)_S\ge C^2\ge 2.$$
Thus $\alpha=2\xi\ge 4$. Theorem \ref{Key} implies that $\varphi_5$ is birational.

Case 2. If $d=2$, then we may take $\beta=1$ whenever $p_g(X)=3$ and $\beta=2$ whenever 
$p_g(X)\ge 4$. Noting that $|G|$ is composed of a pencil, we have 
$G\equiv qC$ where $q\ge p_g(X)-2$. 

When $|G|$ is a rational pencil, since $\pi^*(K_X)|_S
\ge C$, Theorem \ref{Key}(i) is satisfied. 

When $|G|$ is an irrational pencil, then $q\ge 2$. Pick up two different generic 
irreducible  elements $C_1$ and $C_2$ in $|G|$. Then $G-C_1-C_2$ is nef. 
Note that 
$$K_S+\roundup{3\pi^*(K_X)|_S}\ge 
K_S+\roundup{2\pi^*(K_X)|_S}+G.$$
Applying the vanishing theorem, we have the surjective map
$$H^0(S, K_S+\roundup{2\pi^*(K_X)|_S}+G)\rightarrow 
H^0(C_1, K_{C_1}+D_1)\oplus H^0(C_2, K_{C_2}+D_2),$$
where $D_i$ are of positive degree for all $i$. Thus $H^0(C_i, K_{C_i}+D_i)\neq \emptyset$.
So the linear system $|K_S+\roundup{3\pi^*(K_X)|_S}|$ separates different generic irreducible 
elements of $|G|$. 
Theorem \ref{Key}(i) is also satisfied.

By Proposition 5.1 of \cite{Ch2}, $\varphi_4$ is generically finite. This means that 
$|\rounddown{4\pi^*(K_X)}|_S|$ maps a general $C$ onto a curve. Thus $4\pi^*(K_X)|_S\cdot C
\ge 2$ and so $\xi\ge \frac{1}{2}$. 

Now take $m_1=6$. Then $(m_1-3)\xi>1$. Applying Theorem \ref{Key}, one has $\xi\ge \frac{2}{3}$.
Take $m_2=5$. Then $(m_2-3)\xi>1$. Theorem \ref{Key} gives $\xi\ge \frac{4}{5}$. Take $m_3=6$. Then
$(m_3-3)\xi>2$. Theorem \ref{Key} gives $\xi\ge \frac{5}{6}$. Similarly, if we take $m_4=7$, 
then we get $\xi\ge \frac{6}{7}$. 

When $p_g(X)\ge 4$, $\alpha=(5-1-1-\frac{1}{2})\xi>2$. Therefore, by Theorem \ref{Key}, $\varphi_5$
is birational. 

When $p_g(X)=3$, $\alpha=(6-3)\xi>2$. Thus $\varphi_6$ is birational. 
\end{proof}

\begin{rem} Examples 1.3 and 1.5 show that Theorem 3.1 is sharp.
\end{rem}

\begin{thm}\label{T:3.3} Let $X$ be minimal 3-fold of general type with only ${\Bbb Q}$-factorial 
terminal singularities. Assume $p_g(X)\ge 2$, $d=1$ and $g(B)>0$. Then $\varphi_5$ is 
birational onto its image.
\end{thm}
\begin{proof}
This is the simple case. Because $g(B)>0$, the movable part of $|K_X|$ is already basepoint 
free on $X$ and $M_1\equiv aS$ with $a\ge 2$. So one always has $\pi^*(K_X)|_S=\sigma^*(K_{S_0})$, where $S$ is the general fiber 
of the derived fibration $f:X'\rightarrow B$ and $\sigma: S\rightarrow S_0$ is the contraction 
onto minimal model. Note that 
$$\pi^*(K_X)-S-\frac{1}{a}E_1'\equiv (1-\frac{1}{a})\pi^*(K_X).$$
Applying the vanishing theorem, we have the surjective map
$$H^0(X', K_{X'}+\roundup{4\pi^*(K_X)-\frac{1}{a}E_1'})\longrightarrow 
H^0(S, K_S+\roundup{(4-\frac{1}{a})\pi^*(K_X)}|_S).$$
Also note that 
$$K_S+\roundup{(4-\frac{1}{a})\pi^*(K_X)}|_S\ge K_S+3\sigma^*(K_{S_0})+
\roundup{(1-\frac{1}{a})E_1'|_S}.$$

If $(K_{S_0}^2, p_g(S))\neq (1,2)$, then $|K_S+3\sigma^*(K_{S_0})+
\roundup{(1-\frac{1}{a})E_1'|_S}|$ defines a birational map by 2.2 and so does 
$\varphi_5|_S$. Otherwise, because $E_1'|_S\equiv \pi^*(K_X)|_S$ is nef and big and 
effective, we have the same conclusion according to Proposition 2.3.

Because $|5K_{X'}|$ separates different fibers of $f$, we may conclude that $\varphi_5$ 
is birational. 
\end{proof}

\begin{prop}\label{P3.4} Let $X$ be a minimal 3-fold of general type with only 
${\Bbb Q}$-factorial terminal singularities. Assume $p_g(X)\ge 3$, $d=1$ and 
$B={\Bbb P}^1$. Pick up a general fiber $S$ of $f:X'\rightarrow {\Bbb P}^1$ and 
let $\sigma: S\rightarrow S_0$ be the contraction onto minimal model. Suppose that 
there is a movable linear system $|G|$ on $S$ with $G\le 2\sigma^*(K_{S_0})$ and that 
$|G|$ is not composed of an irrational pencil of curves. Then

(i) $|\roundup{4\pi^*(K_X)|_S}|$ separates different generic irreducible elements of $|G|$ 
provided $p_g(X)\ge 4$.

(ii) $|\roundup{5\pi^*(K_X)|_S}|$ separates different generic irreducible elements of $|G|$ 
provided $p_g(X)\ge 3$.
\end{prop}
\begin{proof}
If $p_g(X)\ge 4$, then we have ${\mathcal O}(3)\hookrightarrow f_*\omega_{X'}$. So we get
$${\mathcal O}(1)\otimes f_*\omega_{X'/{\Bbb P}^1}^2\hookrightarrow 
  f_*\omega_{X'}^4.$$
This means that $\varphi_4|_S$ dominates the bicanonical map of $S$. So 
 \newline $\roundup{4\pi^*(K_X)|_S}\ge 2\sigma^*(K_{S_0})$, because $|2\sigma^*(K_{S_0})|$ is basepoint 
free by (2.2.4). Thus $\roundup{4\pi^*(K_X)|_S}\ge G$. We are done.

If $p_g(X)=3$, then we have ${\mathcal O}(2)\hookrightarrow f_*\omega_{X'}$. So we get
$${\mathcal O}(1)\otimes f_*\omega_{X'/{\Bbb P}^1}^2\hookrightarrow 
  f_*\omega_{X'}^5.$$
This means that $\varphi_5|_S$ dominates the bicanonical map of $S$. So \newline 
 $\roundup{5\pi^*(K_X)|_S}\ge 2\sigma^*(K_{S_0})$, because $|2\sigma^*(K_{S_0})|$ is basepoint 
free. Thus $\roundup{5\pi^*(K_X)|_S}\ge G$. We are done.
\end{proof}

\begin{thm}\label{T:3.5} Let $X$  be a minimal 3-fold of general type with only ${\Bbb Q}$-factorial 
terminal singularities. Assume $p_g(X)\ge 3$, $d=1$, $B={\Bbb P}^1$ and  $S$ is a 
surface with invariants $(K_{S_0}^2, p_g(S))\neq (1,2)$ and $(2,3)$ where $S_0$ is the 
minimal model of $S$. Then

(i) $\varphi_5$ is birational whenever $p_g(X)\ge 4$.

(ii) $\varphi_6$ is birational whenever $p_g(X)=3$.
\end{thm}
\begin{proof}
In order to apply Theorem \ref{Key}, we set $G:=2\sigma^*(K_{S_0})$ where 
$\sigma: S\rightarrow S_0$ is the contraction onto minimal model. By (2.2.4), $|G|$ is 
basepoint free. So $|G|$ is not composed of a pencil of curves. A general member 
$C\in |G|$ is actually non-hyperelliptic by (2.2.3).  We set $m=6$, $p=2$ whenever $p_g(X)=3$ 
and $m=5$, $p=3$ whenever $p_g(X)\ge 4$. 

By the previous proposition, Theorem \ref{Key}(i) is satisfied. What we need is a 
suitable value $\beta$ in order to perform the calculation.
Denote by $M_k$ the movable part of $|kK_{X'}|$ for all $k>0$. Note that the total number 
of linear systems we come across in our argument is finite. So we may suppose, by modifying $\pi$, that 
$|M_k|$ is basepoint free for upper bounded $k$. 
Pick up a number $n>0$. Because ${\mathcal O}(p)\hookrightarrow f_*\omega_{X'}$, 
we have 
$${\mathcal O}(1)\otimes f_*\omega_{X'/{{\Bbb P}^1}}^{pn}\hookrightarrow 
f_*\omega_{X'}^{(p+2)n+1}.$$
This means that $M_{(p+2)n+1}|_S\ge pn\sigma^*(K_{S_0})$, because $|pn\sigma^*(K_{S_0})|$
is basepoint free in general. Thus, for all $n>0$, there is an effective ${\Bbb Q}$-divisor
$E_{(n)}$ on $S$ such that
$$((p+2)n+1)\pi^*(K_X)|_S=_{\Bbb Q}(pn)\sigma^*(K_{S_0})+E_{(n)}.$$
So we may choose a $\beta$ in the number sequence $\{\frac{pn}{2(p+2)n+2}|\ n>0\}$. Suppose 
$n$ is sufficient large, one may take a number $\beta\mapsto \frac{p}{2p+4}$, but 
$\beta<\frac{p}{2p+4}$.  

Now if $p_g(X)\ge 4$, then one may choose $\beta$ such that $(5-1-\frac{1}{p}-
\frac{1}{\beta})>0$. By Theorem \ref{Key} and because $C$ is an even divisor, $\varphi_5$ is birational.

If $p_g(X)=3$, we may also choose a $\beta$ such that
$6-1-\frac{1}{2}-\frac{1}{\beta}>0.$
By Theorem \ref{Key}, $\varphi_6$ is birational.
\end{proof}

\begin{thm}\label{T:3.6} Let $X$  be a minimal 3-fold of general type with only ${\Bbb Q}$-factorial 
terminal singularities. Assume $p_g(X)\ge 3$, $d=1$, $B={\Bbb P}^1$ and  $S$ is a 
surface with invariants $(K_{S_0}^2, p_g(S))=(2,3)$ where $S_0$ is the 
minimal model of $S$. Then

(i) $\varphi_5$ is birational whenever $p_g(X)\ge 4$.

(ii) $\varphi_6$ is birational whenever $p_g(X)=3$.
\end{thm}
\begin{proof} We will still apply Theorem \ref{Key}. 
We may take $p=2$ whenever $p_g(X)=3$ and $p=3$ whenever $p_g(X)\ge 4$. Pick up a general 
fiber $S$. Take $|G|$ to be the movable part of $|\sigma^*(K_{S_0})|$. According to 
\cite{BPV}, we know that $|G|$ is not composed of a pencil and a general member $C\in |G|$ 
is a smooth curve of genus 3. By the vanishing theorem, we have the surjective map
$$H^0(X', K_{X'}+\roundup{\pi^*(K_X)-\frac{1}{p}E_1'})\longrightarrow 
H^0(S, K_S+\roundup{\frac{p-1}{p}E_1'}|_S).$$
So we have $2\pi^*(K_X)|_S\ge C$ by Lemma 2.7 of \cite{Ch2}. One may take $\beta=\frac{1}{2}$.  
Note also that $\xi\ge\frac{1}{2}C^2=1$, because $C^2= 2$.

If $p_g(X)\ge 4$, take $m_1=5$ and then $(m_1-1-\frac{1}{p}-\frac{1}{\beta})\xi\ge \frac{5}{3}>1$.
Theorem \ref{Key} gives $\xi\ge \frac{6}{5}$. Take $m_2=6$, we will get $\xi\ge 
\frac{4}{3}$. So $\alpha=(5-1-\frac{1}{p}-\beta)\xi>2$, which means $\varphi_5$ is 
birational.
   
If $p_g(X)=3$, then $\alpha=(6-1-\frac{1}{2}-\beta)\xi>2$. So, by Theorem \ref{Key}, 
$\varphi_6$ is birational.
\end{proof}

\begin{setup}\label{S:3.7} {\bf The last case}. Assume $p_g(X)\ge 3$, $d=1$, 
$B={\Bbb P}^1$ and  $S$ is a surface with invariants $(K_{S_0}^2, p_g(S))=(1,2),$ 
where $S_0$ is the minimal model of $S$. Denote by $\sigma: S\rightarrow S_0$ the 
contraction map. Let $f: X'\rightarrow {\Bbb P}^1$ be the derived fibration. 
Let ${\mathcal L}_0$ be the saturated sub-bundle of $f_*\omega_{X'}$ which is generated 
by $H^0(W, f_*\omega_{X'})$. Because $|K_{X'}|$ is composed of a pencil of surfaces 
and $\varphi_{1}$ factors through $f$, we see that ${\mathcal L}_0$ is a line bundle 
on ${\Bbb P}^1$. Denote ${\mathcal L}_1:=f_*\omega_{X'}/{\mathcal L}_0$. Then we have 
the exact sequence
$$0\longrightarrow  {\mathcal L}_0\longrightarrow  f_*\omega_{X'}
\longrightarrow  {\mathcal L}_1\longrightarrow  0.$$
Noting that $\text{rk}(f_*\omega_{X'})=2$, we see that ${\mathcal L}_1$ is also a line 
bundle. Noting that $H^0({\Bbb P}^1, {\mathcal L}_0)\cong H^0({\Bbb P}^1, f_*\omega_{X'})$,
we have 
$h^1({\Bbb P}^1, {\mathcal L}_0)\ge h^0({\Bbb P}^1, {\mathcal L}_1)$. Note that 
 $\deg({\mathcal L}_0)=p_g(X)-1\ge 2$. We have $h^1({\Bbb P}^1, {\mathcal L}_0)=0$. 
So $h^0({\Bbb P}^1, {\mathcal L}_1)=0$. On the other hand, it's well-known that 
$f_*\omega_{X'/{{\Bbb P}^1}}$ is semi-positive (see \cite{F}). Thus 
$\deg({\mathcal L}_1\otimes\omega_{{\Bbb P}^1}^{-1})\ge 0$. This means 
$\deg({\mathcal L}_1)\ge -2.$ Using the R-R, we may easily deduce that 
$h^1({\mathcal L}_1)\le 1$. So 
$$h^1({\Bbb P}^1, f_*\omega_{X'}))\le 1.$$
So $h^2({\mathcal O}_{X})=h^1({\Bbb P}^1, f_*\omega_{X'})\le 1$. 
\end{setup}

\begin{claim}\label{C:3.8} Keep the same assumption as in 3.7. Fix two 
smooth fibers $S_1$ and $S_2$ of $f$, the restriction map 
$$H^0(X', K_{X'}+S_1+S_2)\longrightarrow H^0(S, K_S)$$
is surjective.
\end{claim}
\begin{proof}
Considering the exact sequence:
$$0\rightarrow {\mathcal O}_{X'}(K_{X'}+S_1-S)\rightarrow 
{\mathcal O}_{X'}(K_{X'}+S_1)\rightarrow {\mathcal O}_{S}(K_S)\rightarrow  0,$$
we have the long exact sequence
\begin{eqnarray*}
\cdots&\rightarrow  H^0(X', K_{X'}+S_1)\overset{\alpha_1}\rightarrow  H^0(S, K_S)
\overset{\beta_1}\rightarrow  H^1(X',K_{X'})\\
&\rightarrow  H^1(X', K_{X'}+S_1)\rightarrow  H^1(S, K_S)=0,
\end{eqnarray*}
If, for a general fiber $S$, $\alpha_1$ is surjective, then
we see that 
$$\dim\Phi_{K_{X'}+S_1}(S)=\dim\Phi_{K_S}(S)=1\ \text{and}\ 
\dim\Phi_{K_{X'}+S_1}(X)=2.$$
So $\dim\Phi_{K_{X'}+S_1+S_2}(X)=2$. We are done. Otherwise, 
$\alpha_1$ is not surjective. Because $\alpha_1\ne 0$, we see that 
$h^2({\mathcal O}_{X'})=h^1(X', K_{X'})\ge 1$. Because $h^2({\mathcal O}_{X'})\le 1$, 
so $h^2({\mathcal O}_{X'})=1$ and $\beta_1$ is surjective. Therefore 
$H^1(X', K_{X'}+S_1)=0$. This also means that $H^1(X', K_{X'}+S')=0$ for any smooth 
fiber $S'$ since $S'\sim S_1$. So we have $H^1(X', K_{X'}+S_1+S_2-S)=0$, which means 
$|K_{X'}+S_1+S_2||_S=|K_S|$. The claim is true.
\end{proof}

\begin{lem}\label{L:3.9} Keep the same assumption as in 3.7. Denote by $|G|$ 
the movable part of $|K_S|$. Pick up a general member $C\in |G|$. Then

(i) $\pi^*(K_X)|_S-\frac{8}{13}C$ is numerically equivalent to an effective 
${\Bbb Q}$-divisor whenever $p_g(X)\ge 4$.

(ii) $\pi^*(K_X)|_S-\frac{7}{12}C$ is numerically equivalent to an effective 
${\Bbb Q}$-divisor whenever $p_g(X)=3$.
\end{lem}
\begin{proof}
We may assume that, on $X'$, all movable parts of a finite number of linear systems are 
basepoint free. Denote by $M_k$ the movable part of $|kK_{X'}|$ for all $k>0$. Denote by 
$M_0$ the movable part of $|K_{X'}+2S_1|$, where $S_1$ is a fixed smooth fiber of $f$. Because 
$\pi^*(K_X)\ge 2S_1$, by  Claim \ref{C:3.8}, we always have 
$$2\pi^*(K_X)|_S\ge M_2|_S\ge M_0|_S\ge G\sim C,$$
where $G$ is the movable part of $|K_S|$ and $C\in |G|$ is a general member. So $C$ is a smooth 
curve of genus two acording to \cite{BPV}.    

Case 1. $p_g(X)\ge 4$.   We consider the sub-system 
$$|K_{X'}+5(K_{X'}+2S_1)+2(K_{X'}+2S_1)+S_1|\subset |13K_{X'}|.$$
Denote by $M_{00}$ the movable part of $|5(K_{X'}+2S_1)|$. By \cite{Ch2}, we know already that 
$\varphi_5$ is generically finite. So $M_{00}$ is nef and big. By the vanishing theorem, we have 
the surjective map
$$H^0(X', K_{X'}+M_{00}+2M_0+S_1)\longrightarrow H^0(S, K_S+M_{00}|_S+2M_0|_S).$$
Noting that $M_{00}|_S\ge 5M_0|_S$, we see that
$M_{13}|_S\ge 8C$ by Lemma 2.7 of \cite{Ch2}, because $|8C|$ is movable on $S$. 
Thus $13\pi^*(K_X)|_S\ge 8C$.

Case 2. $p_g(X)=3$. First, we have the surjective map by the vanishing theorem
$$H^0(X', K_{X'}+\roundup{2\pi^*(K_X)}+S_1)\longrightarrow 
H^0(S, K_S+\roundup{2\pi^*(K_X)}|_S).$$
Denote by $M_{01}$ the movable part of $|3K_{X'}+S_1|$. Then 
we have $M_{01}|_S\ge 2C$ by Lemma 2.7 of \cite{Ch2}. We consider the sub-system
$$|K_{X'}+3(3K_{X'}+S_1)+S_1|\subset |12K_{X'}|.$$
Because $\varphi_9$ is generically finite, the movable part $M_{000}$ of $|3(3K_{X'}+S_1)|$  
is nef and big. Thus we have the surjective map
$$H^0(K_{X'}+M_{000}+S_1)\longrightarrow H^0(S, K_S+M_{000}|_S).$$
Noting that $M_{000}\ge 3M_{01}$, we see that
$12\pi^*(K_X)|_S\ge M_{12}|_S\ge 7C$. We are done.
\end{proof}

\begin{thm}\label{T:3.10} Let $X$  be a minimal 3-fold of general type with only ${\Bbb Q}$-factorial 
terminal singularities. Assume $p_g(X)\ge 3$, $d=1$, $B={\Bbb P}^1$ and  $S$ is a 
surface with invariants $(K_{S_0}^2, p_g(S))=(1,2)$, where $S_0$ is the 
minimal model of $S$. Then

(i) $\varphi_5$ is birational whenever $p_g(X)\ge 4$.

(ii) $\varphi_6$ is birational whenever $p_g(X)=3$.
\end{thm}
\begin{proof}
By Proposition 3.4, Theorem \ref{Key}(i) is satisfied. {}From the claim of
 Proposition 5.3 of \cite{Ch2}, we know that $\xi\ge \frac{3}{5}$. 

If $p_g(X)\ge 4$, we may take $p=3$ and $\beta=\frac{8}{13}$ by Lemma 3.9. 
Take $m_1=5$. Then $(5-1-\frac{1}{p}-\frac{13}{8})\xi\ge \frac{49}{40}>1$. Theorem 
\ref{Key} gives
$\xi\ge \frac{4}{5}$. Take $m_2=6$. Then Theorem \ref{Key} gives $\xi\ge \frac{5}{6}$.
$\cdots$. Take $m_{46}=50$. Theorem \ref{Key} gives $\xi\ge \frac{49}{50}$. 
Now $\alpha=(5-1-\frac{1}{3}-\frac{13}{8})\xi\ge \frac{2431}{1200}>2$. Theorem \ref{Key} 
implies that $\varphi_5$ is birational. 

If $p_g(X)=3$, we may take $p=2$ and $\beta=\frac{7}{12}$ by Lemma 3.9.
Take $m_1=6$. Then $(m_1-1-\frac{1}{2}-\frac{12}{7})\xi\ge \frac{87}{70}$. Theorem 
\ref{Key} gives $\xi\ge \frac{2}{3}$.  
Take $m_2=5$. Then $(5-1-\frac{1}{2}-\frac{12}{7})\xi\ge \frac{25}{21}>1$. 
Theorem \ref{Key} gives $\xi\ge \frac{4}{5}$. Now 
$\alpha=(6-1-\frac{1}{2}-\frac{12}{7})\xi\ge \frac{78}{35}>2$. Thus $\varphi_6$ is 
birational.
\end{proof}

Theorems \ref{T:3.1}, \ref{T:3.3}, \ref{T:3.5}, \ref{T:3.6} and \ref{T:3.10} imply 
Theorem \ref{Main}. 

\section{\bf Applications to the case $p_g=2$}
We present an application of our method to the case $p_g=2$. Throughout this section, 
we always assume that $X$ is a minimal 3-fold of general type with $p_g=2$. Naturally, 
$|K_X|$ is composed of a pencil of surfaces. We keep the same notations as in 2.5. So we have 
a derived fibration $f:X'\rightarrow B$. 

\begin{setup}\label{4.1}  Suppose $g(B)>0$. By Theorem 3.3, $\varphi_5$ is birational. So we 
only need to study the situation with $g(B)=0$.
\end{setup} 

\begin{setup}\label{4.2} {\bf Generalized birationality principle.}  We formulate a generalized 
birationality principle here in order to state the birationality of $\varphi_7$. 
\medskip

\noindent {\bf (GP).} Let $W$ be a smooth projective variety. Suppose that we have a linear system $\Lambda$ which is not 
necessarily a complete one and that we have a movable linear system $|G|$ on $W$. 
Denote by $H$ a generic irreducible element of $|G|$. Assume that

(i) either $H$ is smooth or $H$ is a curve (not necessarily smooth);

(ii) $\Lambda$ separates different generic irreducible elements of $|G|$;

(iii) $\Lambda|_H$ gives a birational map onto its image.

Then $\Lambda$ gives a birational map. 
\medskip

Now we consider the system $|7K_{X'}|$. Let $S$ be a general fiber of $f$. Suppose $g(B)=0$. 
In order to prove the birationality of $\varphi_7$, we have to study the system $|7K_{X'}||_S$
which is not necessarily a complete linear system. On the surface $S$, we always take a 
movable system $|G|$ where either $G:=2\sigma^*(K_{S_0})$ or $G\le 2\sigma^*(K_{S_0})$ and 
$|G|$ is not composed of an irrational pencil of curves. Because 
$${\mathcal O}(1)\otimes f_*\omega_{X'/B}^2\hookrightarrow f_*\omega_{X'}^7,$$
so $|G|\subset |7K_{X'}||_S$. This means that (ii) is satisfied. Note also that (i) is 
automatically satisfied. We may replace Theorem \ref{Key}(i) by (GP) and verify Theorem 
\ref{Key}(ii) and Theorem \ref{Key}(iii) case by case as follows. 
\end{setup}

\begin{setup}\label{4.3} {\bf Finding $\beta$.} We have 
$${\mathcal O}(1)\otimes f_*\omega_{X'/B}^4\hookrightarrow f_*\omega_{X'}^{13}.$$
Because $|4\sigma^*(K_{S_0})|$ is base point free and by Lemma 2.7 of \cite{Ch2}, we have
$M_{13}|_S\ge 4\sigma^*(K_{S_0})$. We still call for a good $\pi$ such that a finite number 
of movable linear systems on $X'$ are all base point free. 
Noting that $M_{13}$ is nef and big, by the vanishing 
theorem, we have the surjective map
$$H^0(X', K_{X'}+M_{13}+S)\longrightarrow  H^0(K_S+M_{13}|_S).$$
Using Lemma 2.7 of \cite{Ch2} again, we have $M_{15}|_S\ge 5\sigma^*(K_{S_0}).$
Repeating this procedure, we may get, for all $t>0$, 
$$(2t+13)\pi^*(K_X)|_S\ge M_{2t+13}|_S\ge (t+4)\sigma^*(K_{S_0}).$$
This means, in fact, that we may choose a $\beta$ such that $\beta\longmapsto \frac{1}{2}$ and
that $\pi^*(K_X)|_S-\beta\sigma^*(K_{S_0})$ is numerically equivalent to an effective 
${\Bbb Q}$-divisor. 
\end{setup}

\begin{setup}\label{4.4} Suppose $(K_{S_0}^2, p_g(S))\neq (1,2)$ and $(2,3)$. Set 
$m=7$, $p=1$ and $G:=2\sigma^*(K_{S_0})$. Then a general $C\in |G|$ is a smooth 
non-hyperelliptic curve. By (4.3), we may take $\beta\mapsto \frac{1}{4}$. So 
$m-1-1-\frac{1}{\beta}>0$. 
Thus Theorem \ref{Key} gives the birationality of $\varphi_7$. 
\end{setup}

\begin{setup}\label{4.5} Suppose $(K_{S_0}^2, p_g(S))=(2,3).$   We set $m=7$, $p=1$ and
$G$ to be the movable part of $|K_S|$. Then $|G|$ is not composed of a pencil and a general 
member $C\in |G|$ is a curve of genus 3 and $C^2=2$. By (4.3), we may take $\beta\mapsto 
\frac{1}{2}$. Also we have $\xi\ge \frac{1}{2}C^2=1$. 
Now $\alpha=(m-1-1-\frac{1}{\beta})\xi>2$. Thus $\varphi_7$ is birational. 
\end{setup}

\begin{setup}\label{4.6} Suppose $(K_{S_0}^2, p_g(S))=(1,2).$  We set $m=8$, $p=1$ and $G$ to be the movable 
part of $|K_S|$. Then a general member $C\in |G|$ is a smooth curve of genus 2. By (4.3), we may take 
$\beta\mapsto \frac{1}{2}$. By the Claim of Proposition 5.3 of \cite{Ch2}, we know $\xi\ge 
\frac{3}{5}$.
Now $\alpha=(m-1-1-\frac{1}{\beta})\xi \mapsto \frac{12}{5}>2$. Thus $\varphi_8$ is birational.
\end{setup}

\begin{setup}\label{4.7}{\bf Summary and remarks.}
What we get is the birationality of $\varphi_8$ for 3-folds of general type with $p_g=2$. 
Unfortunately, we don't know whether this is optimal. We don't have any supporting examples.
\end{setup}

\section*{\bf Acknowledgement}
This note was written while I was visiting the Universit{$\ddot a$}t Essen in the winter of 2002, 
which is subject to the joint German-Chinese project "Komplexe Geometrie" supported by 
the DFG and the NSFC. I am very grateful for the hospitality as well as encouragements from 
both Professor H\'el\`ene Esnault and Professor Eckart Viehweg. This paper benefited a lot 
from my frequent discussions with Professor Viehweg whom I would like to thank for his 
generous helps.  Especially, Professor Viehweg pointed out the fact (2.1.3) to me, which 
greatly simplifies the proof of Theorem 3.5.  
{}Finally I appreciate valuable comments from the referee who points out several unclear arguments 
in the proof of Theorem 2.6 in the original draft.
 


\begin{thebibliography}{99}
\bibitem{BPV} W. Barth, C. Peter, A. Van de Ven, {\em Compact complex surface}, Springer-Verlag, 1984.
\bibitem{Bom} E. Bombieri, {\em Canonical models of surfaces of general type},
Publications I.H.E.S. {\bf 42}(1973), 171-219.
\bibitem{Catan} F. Catanese, {\em Surfaces with $K^2=p_g=1$ and their period mapping},
Springer Lecture Notes in Math. {\bf 732}(1979), 1-29.
\bibitem{Ch1} M. Chen, {\em On pluricanonical maps for threefolds of general type},
J. Math. Soc. Japan {\bf 50}(1998), 615-621. 
\bibitem{Ch2} -----, {\em Canonical stability in terms of singularity index for algebraic 
threefolds}, Math. Proc. Camb. Phil. Soc. {\bf 131}(2001), 241-264.
\bibitem{C-G} S. Chiaruttini, R. Gattazzo, {\em Examples of birationality of pluricanonical 
maps}, Rend. Sem. Mat. Univ. Padova {\bf 107}(2002), 81-94.
\bibitem{Ci} C. Ciliberto, {\em The bicanonical map for surfaces of general type}, 
Proc. Symposia in Pure Math. {\bf 62}(1997), 57-83.
\bibitem{E-L} L. Ein, R. Lazarsfeld, {\em Global generation of pluricanonical
and adjoint linear systems on smooth projective threefolds}, J. Amer. Math. Soc.
{\bf 6}(1993), 875-903.
\bibitem{E-V} H. Esnault, E. Viehweg, {\em Lectures on Vanishing Theorems}. DMV-Seminar 
{\bf 20}(1992), Birkh{$\ddot a$}user, Basel-Boston-Berlin.
\bibitem{F} T. Fujita, {\em On Kahler fiber spaces over curves}, J. Math. Soc. Japan {\bf 30} (1978), 779-794.
\bibitem{Hi} H. Hironaka, {\em Resolution of singularities of an algebraic variety over a field of characteristic zero}, I, Ann. of Math. {\bf 79}(1964), 109-203, II, ibid., 205-326.
\bibitem{Hunt} B. Hunt, {\em Complex manifold geography in dimension 2 and 3}, J. 
Differential Geom. {\bf 30}(1989), 51-153.
\bibitem{Ka} Y. Kawamata, {\em A generalization of Kodaira-Ramanujam's
vanishing theorem}, Math. Ann. {\bf 261}(1982), 43-46.
\bibitem{KMM} Y. Kawamata, K. Matsuda, K. Matsuki, {\em Introduction to the
minimal model problem}, Adv. Stud. Pure Math. {\bf 10}(1987), 283-360.
\bibitem{Kol} J. Koll\'ar, {\em Higher direct images of dualizing sheaves}, I,
Ann. of Math. {\bf 123}(1986), 11-42;\ \ II, ibid. {\bf 124}(1986), 171-202.
\bibitem{K-M} J. Koll\'ar, S. Mori, Birational geometry of algebraic
varieties, 1998, Cambridge Univ. Press.
\bibitem{Lee} S. Lee, {\em Remarks on the pluricanonical and adjoint linear
series on projective threefolds}, Commun. Algebra {\bf 27}(1999), 4459-4476.
\bibitem{Luo} T. Luo, {\em Plurigenera of regular threefolds}, Math. Z. {\bf 217}(1994),
 37-46.
\bibitem{Mi} Y. Miyaoka, {\em The pseudo-effectivity of $3c_2-c_1^2$ for varieties with numerically effective canonical classes}, Algebraic Geometry, Sendai, 1985.
Adv. Stud. Pure Math. {\bf 10}(1987), 449-476.
\bibitem{Sakai} F. Sakai, {\em Weil divisors on normal surfaces}, Duke Math. J.
 {\bf 51}(1984), 877-887.
\bibitem{Sho} V. Shokurov, {\em 3-fold log flips}, Izv. Russ. A. N. Ser. 
Mat. {\bf 56}(1992), 105-203.
\bibitem{V1} E. Viehweg, {\em Vanishing theorems}, J. reine angew. Math. {\bf
335}(1982), 1-8.
\bibitem{X} G. Xiao, {\em Linear bound for abelian automorphisms of varieties
of general type},  J. reine angew. Math. {\bf 476}(1996), 201-207.
\end{thebibliography}
\end{document}